\documentclass[a4paper,11pt]{article}
\usepackage{amsfonts,amssymb,latexsym,amsmath}
\usepackage{multicol}
\usepackage[cp1251]{inputenc}
\usepackage[russian]{babel}
\begin{document}
\Large

\newcommand{\gog}{{\mathfrak g}}
\newcommand{\nog}{{\mathfrak n}}
\newcommand{\mog}{{\mathfrak m}}
\newcommand{\hog}{{\mathfrak h}}
\newcommand{\dog}{{\mathfrak d}}
\newcommand{\bog}{{\mathfrak b}}
\newcommand{\ut}{{\mathfrak u}{\mathfrak t}}
\newcommand{\gl}{{\mathfrak g}{\mathfrak l}}

\newcommand{\GL}{\mathrm{GL}}
\newcommand{\UT}{\mathrm{UT}}
\newcommand{\eps}{\varepsilon}
\newcommand{\ad}{{\mathrm{ad}}}
\newcommand{\Ad}{{\mathrm{Ad}}}
\newcommand{\Ab}{{\Bbb A}}
\newcommand{\BC}{{\cal B}}

\newcommand{\ZC}{{\cal Z}}
\newcommand{\AC}{{\cal A}}
\newcommand{\DC}{{\cal D}}
\newcommand{\YC}{{\cal Y}}
\newcommand{\LC}{{\cal L}}

\newcommand{\BF}{B_i^{(1)}}
\newcommand{\BFF}{B_i^{(1.1)}}
\newcommand{\BFS}{B_i^{(1.2)}}
\newcommand{\BS}{B_i^{(2)}}
\newcommand{\BT}{B_i^{(3)}}
\newcommand{\BFO}{B_i^{(4)}}

\newcommand{\SC}{{\cal S}}
\newcommand{\UC}{{\cal U}}
\newcommand{\FC}{{\cal F}}

\date{}
\title{On index of certain nilpotent Lie algebras
}
\author{Panov A.N.\\ apanov@list.ru
\thanks{The paper is supported by RFFI grants  05-01-00313 and
06-01-00037}}
 \maketitle

\section*{Preliminary information and statements of theorems}

The coadjoint orbits play an important role in representation
theory, simplectic geometry and mathematical physics. According to
the orbit method of A.A.Kirillov ~\cite{K-Orb, K-62} for nilpotent
Lie groups their exists one to one correspondence between the
coadjoint orbits and the irreducible unitary representations in
Hilbert spaces. This gives a possibility to solve the problems of
representation theory and harmonic analysis in geometrical terms of
the orbit space. However the problem of classification of all
coadjoint orbits for specific Lie groups (such as the group of
unitriangular matrices) is an open problem up today that is far from
its solution. In the origin paper on the orbit method ~\cite{K-62}
the description of algebra of invariants and classification of
orbits of maximal dimension was obtained.

In our paper we study the coadjoint representation for nilpotent Lie
algebras that are factors of the  unitriangular Lie algebra by
ideals generated by subsets of root vectors. By this algebra $\LC$
we construct the diagram  $\DC_\LC$ using the formal rule  of
placing of symbols in the table.

Applying this diagram  $\DC$ one can easily calculate the index of
concerned Lie algebra( see theorem 2(3)). Recall that the index of a
Lie algebra is a minimal dimension of stabilizer of a linear form on
this  Lie algebra. For a nilpotent Lie algebra the field of
invariants of the coadjoint representation is a pure transcendental
extension of the main field of degree that is equal to the index of
this algebra ~\cite{Dix}. Respectively, by the diagram one can
easily calculate the maximal dimension of coadjoint orbits (see
theorem 2(2)). We also introduce the method of construction of
system of generators in the field of invariants of the coadjoint
representation (see theorem 1). The main results of the paper is
stated in theorems 1 and 2  in the sequel of this section. .

Note that the author used the diagram in the previous papers for a
classification of all coadjoint orbits for $n\le 7$ ~\cite{P1}, and
also for a description of specific orbits for an arbitrary $n$ (
subragular orbits ~\cite{P1}; orbits associated with involutions
~\cite{P2}).

Let $N=\UT(n,K)$ be a group of lower triangular matrices of size
$n\times n$  with units on the diagonal and with entries in the
field $K$ of zero characteristic. The Lie algebra $\nog=\ut(n,K)$ of
this group consists of all lower triangular matrices of size
$n\times n$ with zeros on the diagonal. The group $N$ acts on the
conjugate space
 $\nog^*$ by the formula
  $\Ad_g^*f(x)=f(\Ad_g^{-1} x)$.
This representation is called  the coadjoint representation.

We identify the symmetric algebra $S(\nog)$ with the algebra of
regular functions   $K[\nog^*]$ on the conjugate space  $\nog^*$. We
also identify $\nog^*$ with the subspace of upper triagular matrices
with zeros on the diagonal. The pairing  of $ \nog$ and  $\nog^*$ is
realized by the Killing form  $(a,b)= \mathrm{Tr}(ab)$, where $a\in
\nog$, $ b\in
 \nog^*$. By this identification the coadjoint
action is realized by the formula $\Ad_g^*b=P(\Ad_g b)$, where $P$
is a natural projection of the space of  $n\times n$-matrices onto
$\nog^*$.

Recall that for any Lie algebra $\gog$ the algebra  $K[\gog^*]$ is a
Poisson algebra with respect to the Poisson bracket such that
 $\{x,y\}=[x,y]$,
~$ x,y\in\gog$. The simplectic leaves of this bracket coincides with
the orbits of coadjoint representation ~\cite{K-Orb}. Respectively,
the algebra of Carimir elements of $K[\gog^*]$ coincides with the
algebra of invariants  $K[\gog^*]^L$ of coadjoint representation.

We also recall that  the coadjoint orbits of an arburary nilpotent
Lie group are the closed subsets with respect to the Zariski
topology in  $\gog^*$ ~\cite[11.2.4]{Dix}.

Consider the standard basis  $\{y_{ij}:~ n\ge i
> j > 1\}$ in the algebra  $\nog$. Denote by $A$ the set of all pairs  $(i,j)$, where
$i>j$. We shall also use the notation  $y_\xi$ for $y_{ij}$, where
$\xi= (i,j)$.

Consider the ideal   $\mog$ in $\nog$  spanned (as a linear
subspace) over the field  $K$ by
 some system of root vectors $\{y_{ij},~ (i,j)\in M\}$, where $M\subset A$.
Denote by $\LC$ the factor algebra  $\nog/\mog$ and by $L$
 the corresponding factor group of $N$ with respect to  the normal subgroup $\exp(\mog)$.

Note that the conjugate space  $\LC^*$ is a subspace of  $\nog^*$
that consists of all  $f\in\nog^*$ which  annihilate on  $\mog$. A
coadjoint  $L$-orbit for  $f\in\LC^*$ coincides with its $N$-orbit.

Consider the order relation  $\succ$ on the set of pairs  $A$ such
that
$$(n,1)\succ (n-1,1)\succ\ldots\succ
 (2,1)\succ(n,2)\succ\ldots\succ(3,2)\succ\ldots\succ(n,n-1).$$

By the ideal  $\mog$ we construct the diagram which is a $n\times
n$-matrix with empty places $(i,j)$, ~$i\le j$ ; the other  places
(i.e. places of  $A$) are filled by symbols  $\otimes$,~$\bullet$,
"$+$" and  "$-$" following the rules stated below. The places
$(i,j)\in M$ are filled by the symbol $\bullet$. We call this
procedure the zero step of the construction of diagram.

Further, we put the symbol $\otimes$ on the greatest place with
respect to the order $\succ$ in  $A\setminus M$. Note that this
symbol with be placed in the first column if  the set of the pairs
of type $(i,1)$ from $A\setminus M$ is nonempty. Assume that we put
the symbol
 $\otimes$ on the place  $(k,t)$,~ $k>t$. Further, we put the symbol
"$-$" on all places  $(k,a)$, ~ $t< a <k$, and  we put the symbol
"$+$" on all places  $(b,t)$,~ $1< b <k$. We finish the first step
of constructing the diagram.

If after this procedure some places of $A$ are not filled we again
put the symbol  $\otimes$ on the greatest (with respect to the order
$\succ$) empty place in $A$. Then,  similarly we  put the symbols
"$+$"\ and "$-$"\ on the empty places taking into account the
following:  we put the symbols  "$+$"\ and "$-$"\ in pairs; if the
both places  $(k,a)$ and  $(a,t)$, where $ k>a>t$ are empty, we put
"$-$"\ on the first place and "$+$"\ on the second place; if one of
these places, $(k,a)$ or $(a,t)$, are  already filled, then we do
not fill the other place. After this procedure we finish the step
that we call a second step.

Continuing the procedure further we have got  the diagram. We denote
this diagram by  $\DC_\LC$. The number of last step is
equal to the number of symbols  $\otimes$ in the diagram.\\
{\bf Example}. Let $n=7$, $\mog = Ky_{51}\oplus Ky_{61}\oplus
Ky_{71}\oplus Ky_{62}$. The corresponding diagram is as follows
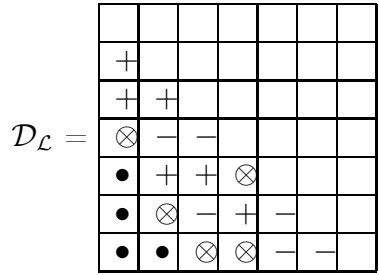
\begin{center}
{\large $\DC_\LC$ =
\begin{tabular}
{|p{0.1cm}|p{0.1cm}|p{0.1cm}|p{0.1cm}|p{0.1cm}|p{0.1cm}|p{0.1cm}|}
\hline  &  &  &  & & &  \\
\hline $+$& & &  & & &  \\
\hline $+$&$+$ & & & & &  \\
\hline $\otimes$ & $-$ & $-$  &  & &  & \\
\hline $\bullet$ & $+$ & $+$&$\otimes$ & & &\\
\hline  $\bullet$ & $\otimes$ & $-$ &$+$  &$-$ & &\\
\hline   $\bullet$&  $\bullet$ & $\otimes$ &$\otimes$  &$-$  &$-$ & \\
\hline
\end{tabular}.}
\end{center}

We construct the diagram in 5 steps, beginning  with zero step:

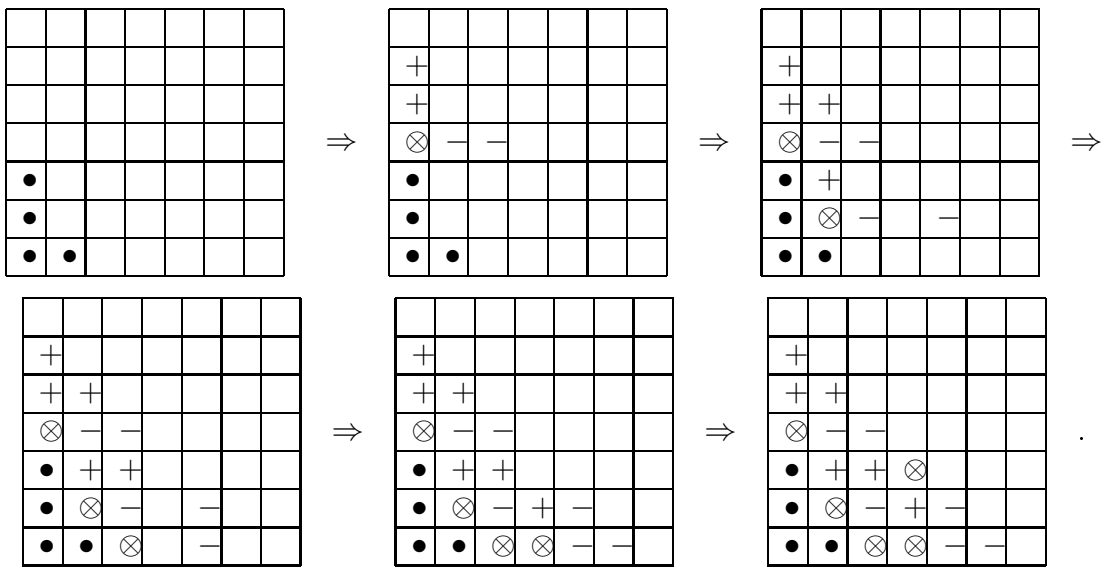
\begin{center}
{\large
\begin{tabular}{|p{0.1cm}|p{0.1cm}|p{0.1cm}|p{0.1cm}|p{0.1cm}|p{0.1cm}|p{0.1cm}|}
\hline  &  &  &  & & &  \\
\hline & & &  & & &  \\
\hline & & & & & &  \\
\hline  &  &   &  & &  & \\
\hline $\bullet$ &  & & & & &\\
\hline  $\bullet$ &  &  &  & & &\\
\hline   $\bullet$&  $\bullet$ &  &  &  & & \\
\hline
\end{tabular}
\quad $\Rightarrow$\quad
\begin{tabular}{|p{0.1cm}|p{0.1cm}|p{0.1cm}|p{0.1cm}|p{0.1cm}|p{0.1cm}|p{0.1cm}|}
\hline  &  &  &  & & &  \\
\hline $+$& & &  & & &  \\
\hline $+$& & & & & &  \\
\hline $\otimes$ & $-$ & $-$  &  & &  & \\
\hline $\bullet$ &  & & & & &\\
\hline  $\bullet$ &  &  &  & & &\\
\hline   $\bullet$&  $\bullet$ &  &  &  & & \\
\hline
\end{tabular}\quad $\Rightarrow$\quad
 \begin{tabular}{|p{0.1cm}|p{0.1cm}|p{0.1cm}|p{0.1cm}|p{0.1cm}|p{0.1cm}|p{0.1cm}|}
\hline  &  &  &  & & &  \\
\hline $+$& & &  & & &  \\
\hline $+$&$+$ & & & & &  \\
\hline $\otimes$ & $-$ & $-$  &  & &  & \\
\hline $\bullet$ & $+$ & & & & &\\
\hline  $\bullet$ & $\otimes$ & $-$ &  &$-$ & &\\
\hline   $\bullet$&  $\bullet$ &  &  &  & & \\
\hline
\end{tabular}\quad $\Rightarrow$}

\end{center}

\begin{center}
{\large
\begin{tabular}{|p{0.1cm}|p{0.1cm}|p{0.1cm}|p{0.1cm}|p{0.1cm}|p{0.1cm}|p{0.1cm}|}
\hline  &  &  &  & & &  \\
\hline $+$& & &  & & &  \\
\hline $+$&$+$ & & & & &  \\
\hline $\otimes$ & $-$ & $-$  &  & &  & \\
\hline $\bullet$ & $+$ & $+$& & & &\\
\hline  $\bullet$ & $\otimes$ & $-$ &  &$-$ & &\\
\hline   $\bullet$&  $\bullet$ & $\otimes$ &  &$-$  & & \\
\hline
\end{tabular}\quad $\Rightarrow$\quad
\begin{tabular}{|p{0.1cm}|p{0.1cm}|p{0.1cm}|p{0.1cm}|p{0.1cm}|p{0.1cm}|p{0.1cm}|}
\hline  &  &  &  & & &  \\
\hline $+$& & &  & & &  \\
\hline $+$&$+$ & & & & &  \\
\hline $\otimes$ & $-$ & $-$  &  & &  & \\
\hline $\bullet$ & $+$ & $+$& & & &\\
\hline  $\bullet$ & $\otimes$ & $-$ &$+$  &$-$ & &\\
\hline   $\bullet$&  $\bullet$ & $\otimes$ &$\otimes$  &$-$  &$-$ & \\
\hline
\end{tabular}\quad $\Rightarrow$\quad
 \begin{tabular}{|p{0.1cm}|p{0.1cm}|p{0.1cm}|p{0.1cm}|p{0.1cm}|p{0.1cm}|p{0.1cm}|}
\hline  &  &  &  & & &  \\
\hline $+$& & &  & & &  \\
\hline $+$&$+$ & & & & &  \\
\hline $\otimes$ & $-$ & $-$  &  & &  & \\
\hline $\bullet$ & $+$ & $+$&$\otimes$ & & &\\
\hline  $\bullet$ & $\otimes$ & $-$ &$+$  &$-$ & &\\
\hline   $\bullet$&  $\bullet$ & $\otimes$ &$\otimes$  &$-$  &$-$ & \\
\hline
\end{tabular}\quad.\quad}
\end{center}

Denote by  $S$ (resp. $C^+$, $C^-$) the set of pairs  $(i,j)$,
filled in the diagram by symbol  $\otimes$ (resp. "$+$", "$-$"). The
set  $A$ of pairs $(i,j)$,~$i>j$ is decomposed into nonintersecting
subsets: $A=M\sqcup C^+\sqcup C^-\sqcup S$.
 We shall give the statement of theorems 1 and 2.

Denote by  $\Ab_m$ the Poisson algebra  $K[p_1,\ldots,p_m;
q_1,\ldots,q_m]$, ~ $\{p_i,q_i\}=1$, and    $\{p_i,q_j\}=0$ for
$i\ne j$.

Recall that a Poisson algebra  $\AC$ decomposes into a product of
two Poisson algebras  $\BC_1\otimes\BC_2$  if  $\AC$ is isomoorphic
as a commutative associative algebra to  $\BC_1\otimes\BC_2$ and
  $\{\BC_1,\BC_2\}=0$.\\
 {\bf Theorem  1}. There exist $z_1,\ldots,z_s\in K[\LC^*]^L$ where
$s=|S|$ such that  \\
1) every  $z_i=y_{\xi_i}Q+P_{>i}$, where $Q$ is some product of
powers of  $z_1,\ldots, z_{i-1}$, and  $P_{>i}$ is a polynomial in
$\{y_\eta\}$,  $\eta\succ \xi_i$;
\\
2) denote by  $\ZC$ the denominator subset, generated by
$z_1,\ldots,z_s$; the localization  $K[\LC^*]_\ZC$ of  algebra
$K[\LC^*]$ with respect to the denominator subset  $\ZC$  is
isomorphic as a Poisson algebra to the tensor product
$K[z_1^\pm,\ldots,z_s^\pm]\otimes \Ab_m$ for some  $m$.

The next section is denoted to the proof of this theorem. We shall
use the procedure of step by step decomposition of Poisson algebras
into tensor products  $\Ab_m\otimes \BC$. Note that in the general
setting not every Poisson algebra  $\AC$, that contains  $\Ab_m$,
can be decomposed into a tensor product of Poisson algebras
$\Ab_m\otimes \BC$. For example, the algebra  $\AC = K[p,q,a]$ with
the bracket $\{p,q\}=1$,~$\{p,a\}=a$ can't be decomposed.
\\
{\bf Theorem  2}. \\
1) The field  $K(\LC^*)^L$ of invariants coincides with the field
$K(z_1,\ldots,z_s)$.\\
2) Maximal dimension of the  coadjoint orbits in $\LC^*$ is equal to
the number of symbols
 "$+$"\ and "$-$"\ in  diagonal $\DC_\LC$.\\
 3) The index of  Lie algebra
 $\LC$ coincides with the number of symbols  $\otimes$ in  diagonal
 $\DC_\LC$.\\
 {\bf Proof of the theorem 2}.  Theorem 1 implies that the algebra
 $K[z_1^\pm,\ldots,z_s^\pm]$ coincides with the algebra of Cazimir elements
 in the Poisson algebra $K[\LC^*]_\ZC$. Hence, $K[z_1^\pm,\ldots,z_s^\pm]$
 is a localization of the algebra of invariants $K[\LC^*]^L$
with respect to  $\ZC$. This follows 1)and further  2) and 3).
$\Box$

\section*{Proof of  Theorem 1}

Denote by  $B_i$ the set of pairs  $(a,b)$,~ $a>b$, unfilled after
the  $i$th step in the procedure of construction of diagram. The
sets  $B_i$ create the chain:
$$B_0\supset B_1\supset\ldots\supset B_s=\emptyset,$$
where $s=|S|$.  Denote $A_i= B_i\sqcup M$,~ $\nog_i=\mathrm{span}\{
y_\eta:~\eta\in A_i\}$,~ $\LC_i=\nog_i/\mog$. Here $A=A_0$.

Let $S=\{\xi_1,\ldots, \xi_s\}$. Recall that  the place  $\xi_i\in
S$ is filled in the diagram by the symbol $\otimes$ during the $i$th
step.

For $1\le i\le s$ we denote by  $C^-_i$ the set of pairs
 $(a,b)$,~$a>b$, that is filled by  "$-$"\ during the  $i$th
step.
\\
{\bf Lemma 1}. The subspace  $\nog_i$ (resp. $\LC_i$) in $\nog$
(resp. $\LC$) is a Lie subalgebra. \\
{\bf Proof}. We shall prove the statement using the method of
induction by $i$. Suppose the the statement is true for all numbers
less than $i$ and the statement is false for number  $i$. Then there
exists $n\ge a>b>c\ge 1 $ such that  $(a,b)\in B_i$, ~~$(b,c)\in
B_i$, and
$$(a,c)\in \bigsqcup_{k\le i} C_k^-\sqcup
\bigsqcup_{k\le i} C_k^+ \sqcup\{\xi_1\ldots \xi_{i}\}.$$

Case $(a,c)\in \{\xi_1\ldots \xi_{i}\}$  is not possible, because if
it occurs then  the places $(a,b)$ and $(b,c)$ are filled by the
symbols "$-$"\ and "$+$"\ respectively; so the both pairs will not
lie
in $B_i$. Consider two last cases.\\
{\bf Case 1.} ~ $(a,c)\in C_k^-$ for some  $k\le i$. Then
$\xi_k=(a,t)$ for some $t<c$. During the  $k$th step the place
$(c,t)$ is filled by the symbol "$+$"\  and  the place $(a,c)$ -- by
the symbol  "$-$".

After the  $(k-1)$th step  the place  $(b,t)$ is filled, otherwise
during the  $k$th step we put the symbol "$+$"\ on the place $(b,t)$
and  the symbol "$-$"\ on the place  $(a,b)$ and obtain $(a,b)\notin
B_i$.

So, after the $(k-1)$th step the places $(c,t)$ and   $(b,c)$ are
empty, and the place  $(b,t)$ is filled by one of the symbols
$\otimes$,~ "$+$"\ or "$-$". This contradicts to the assumption that
$\nog_{k-1}$
is a subalgebra.\\
{\bf Case 2.} ~ $(a,c)\in C_k^+$ for some $k\le i$. This case is
considered similarly. $\Box$

 Denote by $\dog^-_i$, where $1\le i\le s$,  the linear subspace in $\nog$,
 spanned by the vectors
  $y_\eta$ such that  $\xi_i\succ \eta$ и $$\eta
\in\bigsqcup_{1\le j\le i} C^-_j.$$

{\bf Remark}. In the above example:
$$\xi_1=(4,1),\quad \xi_2=(6,2),\quad \xi_3=(7,3), \quad
\xi_4=(7,4),\quad \xi_5=(5,4);$$
$$\begin{array}{ll}
 C^-_1=\{(4,2),(4,3)\},&\quad\quad \dog^-_1 =
 \mathrm{span}\{y_{42},y_{43}\};\\
  C^-_2=\{(6,3),(6,5)\},&\quad\quad  \dog^-_2 =
 \mathrm{span}\{y_{42},y_{43},y_{63}, y_{65}\};\\
  C^-_3=\{(7,5)\}, &\quad\quad \dog^-_3 =
 \mathrm{span}\{y_{43}, y_{63}, y_{65}, y_{75}\};\\
  C^-_4=\{(7,6)\},&\quad\quad \dog^-_4 =
 \mathrm{span}\{y_{65},y_{75},y_{76}\};\\
  C^-_5= \emptyset, &\quad\quad \dog^-_5 = \dog^-_4.\end{array}
  $$
{\bf Lemma 2}. For every $\xi\in S$ the linear subspace $\dog^-_\xi$
is a subalgebra in $\nog$. \\
{\bf Proof}. Suppose that the statement is true for all $l<i$. We
shall prove it for the number $i$. The proof follows from the points
1) and 2) below.
\\
1) Let $(k,b)\in C_i^-$ and $(a,k)\in C_l^-$ for some $l<i$. We
shall show that  $(a,b)$ lies in   $C_m^-$ for some $m<i$.

Assume the contrary. Then  $(a,b)\in B_{i}$ (i.e. the place $(a,b)$
is empty after the  $i$th step). Since $(k,b)\in C_i^-$, then
$\xi_i=(k,t)$ for some $t<b$ and the pair $(b,t)\in C_i^+\subset
B_{i-1}$. So, the pairs  $(a,b)$ and $(b,t)$ lie in $B_{i-1}$. Since
 $\nog_i$ is a subalgebra, then $(a,t)\in M\sqcup B_{i-1}$.

The pair $(a,t)$ can't lie in  $M$ since in this case all
$(a,j)$,~$j\le t$, lie in  $M$. This contradicts    to an existence
of $l<i$ such that $(a,k)\in C_l^-$.

If $(a,t)\in  B_{i-1}$, then $(a,t)\succ (k,t)$ (since $a>k$). We
obtain  a contradiction  to the choice of  $\xi_i=(k,t)$
as a greatest pair in the sense of order relations $\succ$ in $B_{i-1}$.\\
 2) Let $(k,a)\in C_i^-$ and $(a,b)\in C_l^-$ for some  $l<i$.
Let us show that  $(k,b)$ lies in  $C_m^-$ for some  $m<i$.

Assume the contrary. Then $(k,b)\in B_{i}$ (i.e. the place  $(k,b)$
is empty after the  $i$th step). As above  $\xi_i=(k,t)$ for some
$t<b$. The pair  $(b,t)\notin B_{i-1}$. Otherwise during the $i$th
step the place  $(k,b)$ will be filled by the symbol  "$-$", and
$(b,t)$ -- by the symbol "$+$".

So, the place  $(b,t)$ is filled before the  $i$th step.
 According to the  procedure of
placing the symbols, the  symbol $\otimes$ can't take the place
$(b,t)$. The symbol "$-$"\ also can't take the place  $(b,t)$, since
in this case  $(b,t)\in C^-_m$ for some $m<i$ and we have got the
pair of places $(a,b)\in C_l^-$ and $(b,t)\in C^-_m$, where $l,m<i$,
such that  the place $(a,t)$ is empty  after the  $(i-1)$th step
(during the $i$th step we put the symbol "$+$"\ on the place
$(a,t)$). This contradicts to the induction conjecture.

It remains to consider the last case when the place  $(b,t)$ is
filled by the symbol  "$+$". Then there exists a place
$(p,t)$,~$p>k$, that is filled by the symbol $\otimes$ during some
step  $q<i$; at the same time the place  $(p,b)$  is filled by the
symbol "$-$"\ and the places $(p,a)$ and  $(p,k)$ are filled by
"$-$"\ before the  $q$th step. Finally, after the  $(q-1)$th step we
obtain the pair of places  $(p,a)$ and $(a,b)$ filled by the symbol
"$-$" and the empty place  $(p,b)$ (it is filled by the symbol "$-$"
during the  $q$th step). A contradiction. $\Box$

{\bf Proof of the theorem 1}. Consider the localization
$S(\LC_{i-1})_z$ of the algebra  $S(\LC_{i-1})$ with respect to the
denominator subset, generated by  $z=y_{\xi_i}$. Let us show that
for any  $i$ the exists an embedding of Poisson algebras
$$\theta_{i-1}: S(\LC_i)\to S(\LC_{i-1})_z\eqno(1)$$
such that  the Poisson algebra  $S(\LC_{i-1})_z$ decomposes
$$S(\LC_{i-1})_z=\Ab_m\otimes K[z^\pm]\otimes \theta_{i-1}S(\LC_i),\eqno(2)$$
for some  $m$. Then applying the induction by  $i$ we finish the
proof of  2). The proof of 1) follows from the construction of the
embedding $\theta_{i-1}$ that we shall present below.

Let  $\xi_i=(k,t)$,~ $k>t$.  There exists $p>k$ such that the pairs
 of the  $t$th column  from  $M$  fill   $\{(j,t):~
p\le j\le n\}$.

Decompose the set $B_i$ into subsets
$$B_i = \BFF\sqcup\BFS\sqcup\BS\sqcup\BT\sqcup\BFO,\eqno(3)$$
Here every subset is contained of the pairs  $(a,b)\in B_i$ that
satisfy the additional conditions:\\
{\bf 1.1)} ~$(a,b)\in\BFF$, if  $1<a<k$ and $t<b<k$, and  both pairs
 $(a,t)$ and $(k,b)$ lie in $B_{i-1}$;\\
{\bf 1.2)} ~$(a,b)\in\BFS$, if one  these conditions take place:\\
{\bf a)}~ $1<a<k$ and  $b=t$, \\
 {\bf b)}~ $a=k$ and $t<b<k$,\\
 {\bf c)}~ $1<a<k$ and $t<b<k$, at that $(a,t)\notin B_{i-1}$ or
 $(k,b)\notin B_{i-1}$;
 \\
{\bf 2)}~$(a,b)\in\BS$, if $k<a\le n$,~ $t<b<k$;\\
 {\bf 3)} ~$(a,b)\in\BT$, if $p\le a\le n$,~ $b=k$;\\
 {\bf 4)} ~$(a,b)\in\BFO$, if $k<a\le n$,~ $b>k$.
\\
 {\bf Remark}. Note that the pairs  $(a,k)$,~$k<a<p$ do not lie in
 $B_{i-1}$ (precisely,
 they are filled by  the symbol "$-$"\ during the steps with numbers
  $\le i-1$). Note also that all  $(a,j)$, where $p\le a\le n$ and $t<j\le k$,
 lie in $A_{i-1}=B_{i-1}\sqcup M$.

We introduce the following notations:
$$\tilde{y}_{ab}= - \left\vert\begin{array}{cc}y_{at}&y_{ab}\\y_{kt}&y_{kb}
\end{array}\right\vert \cdot y_{kt}^{-1}\quad \mbox{for}\quad
(a,b)\in \BFF;\eqno(4)
$$
$$\tilde{y}_{ak} = \left(y_{ak}y_{kt}+\sum_{(j,t)\in B_{i-1}}
y_{aj}y_{jt}\right)y_{kt}^{-1}\quad \mbox{for}\quad (a,k)\in
\BT;\eqno(5)
$$
$$ \tilde{y}_{ab} = y_{ab}\quad\mbox{для ~всех~ остальных}\quad
(a,b)\in
 B_i.\eqno(6)
$$
We extent the correspondence
$$\theta_{i-1}:~y_{ab}\mapsto \tilde{y}_{ab},\quad (a,b)\in B_i$$
to the embedding (1). By definition,  the subalgebra  $\Ab_m$ is
generated by the elements $p_j=y_{kj}$, ~$(k,j)\in C_i^-$, and
$q_j=y_{jt}y_{kt}^{-1}$, ~$(j,t)\in C_i^+$.

Easy to see that  $S(\LC_{i-1})_z$ decomposes as a commutative
associative algebra  into  the tensor product  $\Ab_m\otimes
K[z^\pm]\otimes \theta_{i-1}S(\LC_i)$. The subalgebras $\Ab_m$,~
$K[z^\pm]$ and $ \theta_{i-1}S(\LC_i)$ are pairwise in
involution(recall that  two elements $u$ and $v$ of  a Poisson
algebra are in involution if $\{u,v\}=0$).

 It remains to show that  $\theta_{i-1}$ is an embedding of Poisson
algebras. It is necessary and sufficient to check that
$$\{\tilde{y}_{ad},\tilde{y}_{bc}\} = \delta_{db}\tilde{y}_{ac}\bmod\mog,\eqno(7)$$
for all  $a>d$,~$b>c$,~$a>c$ where   $(a,d)$ and $(b,c)$ lie in
$B_i$.

We prove the condition (7) by running over all case of entering of
pairs
$(a,d)$, ~$(b,c)$ in the subsets of decomposition  (3).\\
{\bf Case 1}. $(a,d)\in \BFO$. One can prove (7) easily.
\\
{\bf Case  2}.~$(a,d)\in \BT$, ~$(b,c)\in \BS$.\\
In this case,  $d=k$,~ $a\ge p$ and $\{y_{ak},y_{bc}\}=0$. From  (6)
we obtain $\tilde{y}_{bc}=y_{bc}$. We have to show that
$$\{\tilde{y}_{ak},y_{bc}\}=0\bmod\mog.\eqno(8)$$

The element  $y_{bc}$ is in involution with all  elements of the
standard basis  from  (5) apart from $y_{ct}$(if this element really
appears in (5)). If $b\ge p$, then
$\{y_{bc},y_{ct}\}=y_{bt}\in\mog$, this proves  (8).

Let $k<b<p$. Let us show that an existence of  $(b,c)\in B_i$, where
$k<b<p$, implies that  $(c,t)\notin B_{i-1}$. Really, let $(c,t)\in
B_{i-1}$. The pair  $(b,c)$ lies in $B_i$ and, therefore, also it
lies in  $B_{i-1}$. Since  $\nog_i$ is a subalgebra, then $(b,t)\in
B_{i-1}$. On the other side, all places $(j,t)$,~$k<j<p$, are
already filled after $(i-1)$th step. That is no one of these pairs
are contained in  $B_{i-1}$ (see the above remark). A contradiction.
So $(c,t)\notin B_{i-1}$ and, therefore, $y_{ct}$ does not appear in
(5).
The equality  (8) is true.\\
{\bf Case  3}.~$(a,d)\in \BT$, ~$(b,c)\in \BFF$.\\
As above  $d=k$, ~$a\ge p$ and $\{y_{ak},y_{bc}\}=0$.  We have to
prove that
$$\{\tilde{y}_{ak},\tilde{y}_{bc}\}=0\bmod\mog,\eqno(9)$$
where $\tilde{y}_{ak}$ and  $\tilde{y}_{bc}$ from  (5) and  (4)
respectively.

The element  $y_{bc}$  is in involution with all elements of the
standard basis from (5) apart from  $y_{ab}$ and $ y_{ct}$ (if these
elements really appear in (5)). Since  $(b,c)\in \BFF$, then
$y_{bt}\in B_{i-1}$ and the element  $y_{ab}$ really appears in the
sum of  (5).

As for  the element   $ y_{ct}$  we can't say this exactly. If
$(c,t)\in B_{i-1}$, then the element  $ y_{ct}$ appears in the sum
of (5). The equality  (9) is fulfilled since all elements of the
standard basis from (4) are in involution with  $\tilde{y}_{ak}$.

If $(c,t)\notin B_{i-1}$, then the element  $ y_{ct}$ does not
appear in the sum of (5). Since $a\ge p$, then $y_{at}\in\mog$; the
equality  (9) is checked directly:
$$\{\tilde{y}_{ak},\tilde{y}_{bc}\}= -\{y_{ak}y_{kt}+y_{ab}y_{bt},
 \left\vert\begin{array}{cc}y_{bt}&y_{bc}\\y_{kt}&y_{kc}
\end{array}\right\vert \}\cdot y_{kt}^{-2} =  0\bmod\mog.$$
 {\bf Case  4}.~ $(a,d)\in \BT$, ~$(b,c)\in \BFS$.\\
 As above  $d=k$. Recall that  $\BFS$  decomposes into three subsets
 (see. a),~b),~c)
 from  the definition of $\BFS$).\\
{\bf 4a}. $b=k$ and $t<c<k$. Then $\{y_{ak},y_{kc}\}=y_{ac}$. The
element  $y_{kc}$ are in involution with all elements of the
standard basis from (5) apart from  $y_{ak}$ and  $y_{ct}$ (if the
last element really appears in (5)). In our case the place  $(c,t)$
is filled after the  $(i-1)$th step (othewise during the $i$th step
the place  $(c,t)$ will obtain the symbol "$+$"\ and the place
$(k,c)$ -- the symbol "$-$"). The element
$y_{ct}$ does not appear in (5). Hence $\{\tilde{y}_{ak}, y_{kc}\}=0$. \\
{\bf 4b}. $c=t$ and  $1<b<k$. Then  $\{y_{ak},y_{bt}\}=0$. The
element $y_{bt}$ is in involution with all elements of the standard
basis  from  (5) apart from  $y_{ab}$. Since  $y_{at}\in\mog$, then
$$\{\tilde{y}_{ak}, y_{bt}\}=-\{y_{ab}y_{bt},y_{bt}\}y_{kt}^{-1} =
-y_{at}y_{bt}y_{kt}^{-1} = 0\bmod\mog.$$ {\bf 4c}. $1<b<k$ and
$t<c<k$, at that  $(b,t)\notin B_{i-1}$ or $(k,c)\notin B_{i-1}$.
Since in this case  $\{y_{ak},y_{bc}\}=0$ and
$\tilde{y}_{bc}=y_{bc}$, then we have to show that
$$\{\tilde{y}_{ak}, y_{bc}\}=0.\eqno(10)$$

If $(b,t)$ and $(c,t)$ lie in $B_{i-1}$, then the calculation
$$ \{\tilde{y}_{ak}, y_{bc}\} = - \{y_{ab}y_{bt}+y_{ac}y_{ct},
y_{bc}\}y_{kt}^{-1}=0$$ proves (10).

 If  the both pairs  $(b,t)$ and $(c,t)$ do not lie in  $B_{i-1}$,
then $y_{bc}$ is in involution with all elements of the standard
basis from (5). This proves  (10).

The case   $(b,t)\notin B_{i-1}$,~ $(c,t)\in B_{i-1}$ is not
possible, since  in this case  the both pairs  $(b,c)$ and $(c,t)$
lie  $B_{i-1}$, but $(b,t)$ do not. This contradicts to the  fact
that  $\nog_{i-1}$ is a subalgebra.

Let us show that the case  $(b,t)\in B_{i-1}$,~ $(c,t)\notin
B_{i-1}$ is also not possible. Really,  $(b,t)\in B_{i-1}$ and
$(b,c)\in \BFS$ imply    $(k,c)\notin B_{i-1}$. That is the place
$(k,c)$ is filled after the  $(i-1)$th step ( by the symbol "$-$").
Since $(c,t)\notin B_{i-1}$, then the place $(c,t)$ is also filled
after the  $(i-1)$th step. The symbol $\otimes$ can't take this
place because of the procedure of placing the symbols. The symbol
"$-$"\ also can't take this place because in this case after the
$(i-1)$th step we have got symbol "$-$"\ on the places $(k,c)$ and
$(c,t)$, and the place $(k,t)$ is empty(during the $i$th step  this
place is filled by $\otimes$). This contradicts to lemma 2.

It remains  the last case: the place  $(c,t)$ is filled by the
symbol "$+$". There exists a pair $\xi_l =(q,t)\in S$, ~$q>k$,
~$l<i$. According to the procedure of placing of symbols the place
$(q,k)$ is filled by the symbol "$-$"\  before the  $l$th step
(otherwise the place  $(k,t)$ is filled by the symbol  "$+$"\ and
the place $(q,k)$, respectively, by  "$-$"). The symbol "$-$"\
appears on the place  $(k,c)$ also before the $l$th step. We have
got that after the $(l-1)$th step the symbol "$-$" on the places
 $(q,k)$ and  $(k,c)$, but the place $(q,t)$ is empty. This contradicts
 to lemma 2 for
  $i=l-1$.

The equality  (7) is obvious for cases  $(a,d)\in \BS$,
$(b,c)\in\BS$.  Denote  $\BF=\BFF\sqcup\BFS $.\\
{\bf Case 5}. ~$(a,d)\in \BS$, ~$(b,c)\in \BF$.\\
Recall that  here  $\tilde{y}_{ad} = y_{ad}$. The equality (7) if
obvious  for the case   $(b,c)\in\BFS$. Consider the case
$(b,c)\in\BFF$. One can easily check  (7) for the case $d\ne b$ and
for the case $a>p$.

Let us prove that the last case   $d=b$,~ $k<a<p$,~ $(a,b)\in \BS$,~
$(b,c)\in\BFF$ is not possible.

From the above remark we see that all pairs of the form  $(j,k)$,
~$k<j<p$, do not lie in  $B_{i-1}$.  Precisely, all $(j,k)$,
~$k<j<p$, are filled by the symbol  "$-$"\ before the  $(i-1)$th
step. This concerns the  pair  $(a,k)$.

On the other hand, since  $(b,c)\in \BFF$, then after the $i$th step
the place $(b,t)$  is either empty, or is filled by the symbol
"$+$". In any case after the $i$th step the place  $(k,b)$ is filled
by the symbol "$-$".

Finally, after the  $i$th step the places  $(a,k)$ and $(k,b)$ are
filled by the symbol  "$-$".  Then  by  lemma  2  the place $(a,b)$
must be filled by the symbol "$-$"\ during the $(i-1)$th step.
This contradicts to  $(a,b)\in B_i$.\\
{\bf Case 6}.~$(a,d)\in \BF$, ~$(b,c)\in \BF$.\\
The equality  (7) is obvious when two pairs lie in  $\BFS$ and can
be easily checkable in the case when two pairs lie in  $\BFF$. On
can also easily check  (7) in the case $d\ne b$ and in the case when
one of the pairs lies in $\BFS$ and  satisfies the points  a) and
b)(see the definition of $\BFS$).
Consider the last cases. \\
i) Let  $(a,b)\in \BFF$ and $(b,c)\in \BFS$ where $(k,c)\notin
B_{i-1}$ or $(b,t)\notin B_{i-1}$ (see condition 1.2(с)). Let us
show that in this case  $(a,c)\in\BFF$. Then (7) is easily
checkable.

 Suppose the contrary $(a,c)\in \BFS$.
Then one of the places  $(a,t)$ or $(k,c)$ is filled after the
$((i-1))$th step. Since  $(a,b)\in \BFF$ , then the place  $(a,t)$
is empty.
 Hence  the place $(k,c)$ is filled after the $(i-1)$th step. On the other hand,
  $(a,b)\in
\BFF$ implies that  $(k,b)$ is empty  after the  $(i-1)$th step. The
place $(b,c)$ is also empty. Since $\nog_{i-1}$ is a subalgebra,
then  the place  $(k,c)$ is empty after the  $(i-1)$th step. A
contradiction.
\\
ii) $(a,b)\in \BFS$ and $(b,c)\in \BFF$, where $(k,b)\notin B_{i-1}$
or $(a,t)\notin B_{i-1}$ (see condition 1.2(с)).
 Let us show that $(a,c)\in\BFF$ (this implies  (7)).

Really,  $(b,c)\in\BFF$ implies  $(k,c)\in B_{i-1}$. The place
$(a,t)$ can't be filled after the  $(i-1)$th step because at the
same time  the places  $(a,b)$ and  $(b,t)$ are empty. Since  the
pairs  $(a,t)$ and  $(k,c)$ lie in $B_{i-1}$ we have got $(a,c)\in
\BFF$.
 $\Box$

\end{document}